\documentclass[12pt]{amsart}
\usepackage[left=2.2cm,top=3.3cm,right=2.2cm]{geometry}
\usepackage{amsmath}
\usepackage{amssymb}
\usepackage{setspace}
\usepackage{enumerate}
\begin{document}

\newtheorem*{EVT}{Extreme Value Theorem}

\title{A One-Sentence Line-of-Sight Proof of the Extreme Value Theorem}
\author{Samuel J. Ferguson}
\address{Courant Institute of Mathematical Sciences, New York University, 251 Mercer St,
New York, NY 10012}
\email{ferguson@cims.nyu.edu}

\begin{abstract}
We give a one-sentence proof that a continuous real-valued function f on a closed, bounded interval attains a maximum value, by the following device. We define x in [a, b] to be a lookout point if f(t) does not exceed f(x) whenever t lies in [a, x). Letting c be the maximum of the set of lookout points, we prove that f(c) is the maximum value of f.
\end{abstract}

\maketitle

The maximum value of a continuous real-valued function $f$ on $[a, b]$ is attained at its largest ``lookout point." We call $x$ in $[a, b]$ a \emph{lookout point} if, whenever $t$ lies in $[a, x)$, we have
\[
f(t)\leq f(x).
\]
The set $L$ of lookout points is closed. Indeed, let $x_n\to x$, with $x_n$ in $L$. If $t$ is in $[a, x)$, then eventually $t$ lies in $[a, x_n)$, so
\[
f(t)\leq f(x_n).
\]
By continuity, $f(t)\leq f(x)$, as desired. We use the fact that a closed, bounded, and nonempty set has a maximum and a minimum. In particular, $\max(L)$ exists.

\begin{EVT}
If $f$ is a real-valued continuous function on $[a, b]$ then $f$ has a maximum value on $[a, b]$. In other words, for some $c$ in $[a, b]$, no value of $f$ exceeds $f(c)$.
\end{EVT}

\begin{proof}
Letting
\[
L = \{x\textrm{ in }[a, b]\textrm{ such that }t\textrm{ in }[a, x)\textrm{ implies } f(t)\leq f(x)\}
\]
and $c = \max(L)$, it suffices to show that, given $k > f(c)$, the closed, bounded set
\[
S_k = \{t\textrm{ in }[a, b]\textrm{ such that }f(t)\geq k\}
\]
is empty, which follows since, if some $d$ satisfies $f(d)\geq k$, then $d > c$, whence $d$ is not in $L$, so there exists a $t < d$ for which
\[
f(t) > f(d)\geq k,
\]
proving that $S_k$ has no minimum.
\end{proof}












\end{document}